\documentclass[twocolumn]{autart}    

\usepackage{natbib}
\usepackage{graphicx}          

\usepackage{subdepth} 

\usepackage{graphics} 
\usepackage{epsfig} 

\usepackage{amssymb} 
\usepackage{amsmath,psfrag,color,graphicx} 
\usepackage{graphicx}
\usepackage{quoting} 
\usepackage{lettrine} 
\usepackage{tocloft} 
\usepackage{titletoc}
\usepackage{fmtcount}
\usepackage{calc} 
\usepackage{centernot}
\usepackage{nicefrac}
\usepackage{dcolumn}
\usepackage{dsfont}

\usepackage{lipsum}
\usepackage{psfrag}
\usepackage[export]{adjustbox} 

\PassOptionsToPackage{cmyk}{xcolor}
\usepackage[cmyk]{xcolor}
\selectcolormodel{cmyk}

\usepackage{textcomp}

\usepackage{txfonts}
\usepackage{lmodern}
\usepackage{times}

\makeatletter
\let\save@mathaccent\mathaccent
\newcommand*\if@single[3]{%
	\setbox0\hbox{${\mathaccent"0362{#1}}^H$}%
	\setbox2\hbox{${\mathaccent"0362{\kern0pt#1}}^H$}%
	\ifdim\ht0=\ht2 #3\else #2\fi
}
\newcommand*\rel@kern[1]{\kern#1\dimexpr\macc@kerna}
\newcommand*\widebar[1]{\@ifnextchar^{{\wide@bar{#1}{0}}}{\wide@bar{#1}{1}}}
\newcommand*\wide@bar[2]{\if@single{#1}{\wide@bar@{#1}{#2}{1}}{\wide@bar@{#1}{#2}{2}}}
\newcommand*\wide@bar@[3]{%
	\begingroup
	\def\mathaccent##1##2{%
		\let\mathaccent\save@mathaccent
		\if#32 \let\macc@nucleus\first@char \fi
		\setbox\z@\hbox{$\macc@style{\macc@nucleus}_{}$}%
		\setbox\tw@\hbox{$\macc@style{\macc@nucleus}{}_{}$}%
		\dimen@\wd\tw@
		\advance\dimen@-\wd\z@
		\divide\dimen@ 3
		\@tempdima\wd\tw@
		\advance\@tempdima-\scriptspace
		\divide\@tempdima 10
		\advance\dimen@-\@tempdima
		\ifdim\dimen@>\z@ \dimen@0pt\fi
		\rel@kern{0.6}\kern-\dimen@
		\if#31
		\overline{\rel@kern{-0.6}\kern\dimen@\macc@nucleus\rel@kern{0.4}\kern\dimen@}%
		\advance\dimen@0.4\dimexpr\macc@kerna
		\let\final@kern#2%
		\ifdim\dimen@<\z@ \let\final@kern1\fi
		\if\final@kern1 \kern-\dimen@\fi
		\else
		\overline{\rel@kern{-0.6}\kern\dimen@#1}%
		\fi
	}%
	\macc@depth\@ne
	\let\math@bgroup\@empty \let\math@egroup\macc@set@skewchar
	\mathsurround\z@ \frozen@everymath{\mathgroup\macc@group\relax}%
	\macc@set@skewchar\relax
	\let\mathaccentV\macc@nested@a
	\if#31
	\macc@nested@a\relax111{#1}%
	\else
	\def\gobble@till@marker##1\endmarker{}%
	\futurelet\first@char\gobble@till@marker#1\endmarker
	\ifcat\noexpand\first@char A\else
	\def\first@char{}%
	\fi
	\macc@nested@a\relax111{\first@char}%
	\fi
	\endgroup
}
\makeatother

\newcommand{\ve}[2][]{\ensuremath{\boldsymbol{\mathrm{#2}}}_{#1}}
\newcommand{\vet}[2][]{\ensuremath{\smash{\boldsymbol{\mathrm{#2}}^{\!\top}_{#1}}}}
\newcommand{\vd}[2][]{\ensuremath{\dot{\boldsymbol{\mathrm{#2}}}_{#1}}}

\newcommand{\ma}[2][]{\ensuremath{\boldsymbol{\mathrm{#2}}}_{#1}}
\newcommand{\mat}[2][]{\ensuremath{\boldsymbol{\mathrm{#2}}^{\!\top}_{#1}}}
\newcommand{\md}[2][]{\ensuremath{\dot{\boldsymbol{\mathrm{#2}}}}_{#1}}

\newcommand{\maw}[2][]{\ensuremath{\overline{\boldsymbol{\mathrm{#2}}} _{#1}}}

\DeclareMathOperator{\vect}{\ensuremath{\mathrm{vec}}}

\DeclareMathOperator{\skews}{\ensuremath{\mathrm{skew}}}

\DeclareMathOperator{\syms}{\ensuremath{\mathrm{sym}}}

\newcommand{\R}{\ensuremath{\mathds{R}}}

\newcommand{\N}{\ensuremath{\mathds{N}}}

\newcommand{\SOT}{\ensuremath{\mathsf{SO}(3)}}

\newcommand{\SO}{\ensuremath{\mathsf{SO}(n)}}

\newcommand{\so}{\ensuremath{\mathsf{so}(n)}}

\newcommand{\St}{\ensuremath{\mathsf{St}(p,n)}}

\newcommand{\StN}{\ensuremath{(\mathsf{St}(p,n))^N}}

\newcommand{\Sn}{\ensuremath{\mathsf{S}^{n}}}

\newcommand{\GC}{\ensuremath{\mathsf{S}^1}}

\newcommand{\twosphere}{\ensuremath{\mathsf{S}^2}}

\newcommand{\ie}{\textit{i.e.}, }

\newcommand{\eg}{\textit{e.g.}, }

\newcommand{\mtr}{\hspace{-0.3mm}\ensuremath{^\top}} 

\newcommand\raiseT[2]{\setbox0\hbox{$#1{#2}$}\raise\dp0\box0}

\newcommand{\V}{\ensuremath{\mathcal{V}}}

\newcommand{\E}{\ensuremath{\mathcal{E}}}

\newcommand{\etc}{\emph{etc.\ }}

\DeclareMathOperator{\Log}{\ensuremath{\mathrm{Log}}}

\DeclareMathOperator{\trace}{\ensuremath{\mathrm{tr}}}

\DeclareMathOperator{\im}{\ensuremath{\mathrm{Im}}}

\DeclareMathOperator{\re}{\ensuremath{\mathrm{Re}}}

\newcommand{\Ni}{\ensuremath{\mathcal{N}_i}}

\makeatletter
\newcommand{\raisemath}[1]{\mathpalette{\raisem@th{#1}}}
\newcommand{\raisem@th}[3]{\raisebox{#1}{$#2#3$}}
\makeatother

\newcommand{\ts}[2][]{\ensuremath{\mathsf{T}_{#2}#1}}

\newcommand{\AGAS}{\textsc{agas}}

\definecolor{kthbluergb}{RGB}{25,84,166}
\definecolor{kthbluecmyk}{cmyk}{1,0.55,0,0}

\definecolor{kthblueA}{RGB}{25,84,166}
\definecolor{kthblueB}{RGB}{46,124,192}
\definecolor{kthblueC}{RGB}{112,153,209}
\definecolor{kthblueD}{RGB}{164,186,225}
\definecolor{kthblueE}{RGB}{211,220,241}

\newcounter{counter} 

\newtheorem{theorem}[counter]{Theorem}

\newtheorem{proposition}[counter]{Proposition}

\newtheorem{definition}[counter]{Definition}

\newtheorem{remark}[counter]{Remark}


\newcounter{parentnumber}

\allowdisplaybreaks

\begin{document}

\begin{frontmatter}

\title{High-dimensional Kuramoto models on Stiefel manifolds\\ synchronize complex networks  almost globally\thanksref{footnoteinfo}} 

\thanks[footnoteinfo]{Some preliminary results of this paper appear in \cite{markdahl2018cdc}. Johan Markdahl is corresponding author. The work of Johan Markdahl is supported by the University of Luxembourg internal research project \textsc{pppd}.}

\author[XXX]{Johan Markdahl}\ead{markdahl@kth.se},    
\author[YYY]{Johan Thunberg}\ead{johan.thunberg@hh.se},               
\author[XXX]{Jorge Goncalves}\ead{jorge.goncalves@uni.lu}  

\address[XXX]{Luxembourg Centre for Systems Biomedicine, University of Luxembourg, Belval, Luxembourg}  
\address[YYY]{School of Information Technology, Halmstad University, Halmstad, Sweden}

\begin{keyword}                           
Synchronization; Kuramoto model; Stiefel manifold; Multi-agent system; decentralization; networked robotics.               
\end{keyword}                             

\begin{abstract} 
	
\noindent The Kuramoto model of coupled phase oscillators is often used to describe synchronization phenomena in nature. Some applications, \eg  quantum synchronization and rigid-body attitude synchronization, involve high-dimensional Kuramoto models where each oscillator lives on the $n$-sphere or $\SO$. These manifolds are special cases of the compact, real Stiefel manifold $\St$. Using tools from optimization and control theory, we prove that the generalized Kuramoto model on $\St$ converges to a synchronized state for any connected graph and from almost all initial conditions provided $(p,n)$ satisfies $p\leq \tfrac23n-1$ and all oscillator frequencies are equal. This result could not have been predicted based on knowledge of the Kuramoto model in complex networks over the circle. In that case, almost global synchronization is graph dependent; it applies if the network is acyclic or  sufficiently dense. This paper hence identifies a property that distinguishes many high-dimensional generalizations of the Kuramoto models from the original model.	
\end{abstract}

\end{frontmatter}

\section{Introduction}

The Kuramoto model and its many variations are canonical models of systems of coupled phase oscillators \citep{hoppensteadt2012weakly}. As such, they are abstract models that capture the essential properties observed in a wide range of synchronization phenomena. However, many properties of a particular system are lost through the use of these  models. In this paper we study the convergence of a multi-agent system on the Stiefel manifold that includes the Kuramoto model as a special case. For a system of $N$ coupled agents that are subject to various constraints, a high-dimensional Stiefel manifold may provide a more faithful approximation of reality than a phase oscillator model. The orientation of an agent in a swarm can \eg be modeled as an element of the circle, the sphere, or the rotation group---all of which are Stiefel manifolds. For a high-dimensional model to be preferable it must retain some property of the original system which is lost in phase oscillator models. That is indeed the case; we prove that if the complex network of interactions is connected, if all frequencies are equal, and a condition on the parameters of the manifold is satisfied, then the system converges to the set of synchronized states from almost all initial conditions. The same cannot be said about the Kuramoto model in complex networks on the circle $\smash{\mathsf{S}^1}$ in the case of oscillators with homogeneous frequencies \citep{rodrigues2016kuramoto}. Under that model, guaranteed almost global synchronization  requires that the complex network can be represented by a graph that is acyclic or sufficiently dense \citep{dorfler2014synchronization}. To characterize all such graphs is an open problem. 

Since the Stiefel manifold includes the $n$-sphere and the special orthogonal group as special cases, there is a considerable literature on synchronization on particular instances of the Stiefel manifold. Previous works that address synchronization on all Stiefel manifolds is limited to  \cite{thunberg2017dynamic} which relies on the so-called dynamic consensus approach (see \cite{scardovi2007synchronization,sarlette2009consensus}). The dynamic consensus approach is used to stabilize the consensus manifold on $\St$ almost globally for any quasi-strongly connected digraph. However, dynamic consensus requires the introduction of auxiliary variables that are communicated in a second, undirected graph. The gradient descent flow studied in this paper is preferable to \cite{thunberg2017dynamic} in the case of $p\leq\tfrac23n-1$ since it provides the same convergence guarantees but uses less communication and computation. If $p>\tfrac23n-1$, then \cite{thunberg2017dynamic} is preferable. Note that for modeling synchronization in nature the gradient descent flow is arguably always preferable since the auxiliary variables in \cite{thunberg2017dynamic} do not have a physical interpretation.

The problem of almost global synchronization of multi-agent systems on nonlinear spaces has received some attention in the literature, see the survey \cite{sepulchre2011consensus}. Until recently, there have been three main approaches: potential shaping which is based on gradient descent flows \citep{tron2012intrinsic}, probabilistic gossip algorithms \citep{mazzarella2014consensus}, and dynamic consensus algorithms. \cite{markdahl2018tac} shows that a fourth approach based on gradient descent flows, which can be interpreted as high-dimensional Kuramoto models, yields almost global synchronization on the $n$-sphere for all $n\geq2$. It requires less communication and computation, but is limited to undirected graphs and certain manifolds. This paper establishes that it works not just on $\mathcal{S}^n$ but also on $\St$ when $p\leq\tfrac23n-1$.



The Kuramoto model on the $n$-sphere is known as the Lohe model \citep{lohe2010quantum}. 
%
%
Many works on the Lohe model concern the complete graph case \citep{olfati2006swarms,lohe2010quantum,li2014unified,lohe2018higher}. Almost global stability of the consensus manifold in the case of a complete graph and homogeneous frequencies has been shown for the Kuramoto model \citep{watanabe1994constants}, Lohe model \citep{olfati2006swarms}, and on rather general manifolds  \citep{sarlette2009consensus}. The Kuramoto model on networks is less well-behaved \citep{canale2015exotic}. Most results for the Lohe model on networks show convergence from a hemisphere \citep{zhu2013synchronization,thunberg2018lifting,zhang2018equilibria}. Many papers address the case of heterogeneous frequencies  \citep{chi2014emergent,chandra2018continuous,ha2018relaxation}. Some concern the thermodynamic limit $N\rightarrow\infty$, where $N$ denotes the number of agents   
\citep{chi2014emergent,tanaka2014solvable,ha2018relaxation,frouvelle2019long}. There is also a discrete-time model \citep{li2015collective}.

Applications for synchronization on $\smash{\mathcal{S}^2}$ include synchronization of interacting tops \citep{ritort1998solvable}, modeling of collective motion in flocks \citep{al2018gradient}, autonomous reduced attitude synchronization and balancing \citep{song2017intrinsic}, synchronization in planetary scale sensor networks \citep{paley2009stabilization}, and consensus in opinion dynamics \citep{aydogdu2017opinion}. Applications on $\smash{\mathcal{S}^3}$ include synchronization of quantum bits \citep{lohe2010quantum} and models of learning \citep{crnkic2018swarms}. The Kuramoto model on $\SOT$ is of interest in rigid-body attitude synchronization \citep{sarlette2009consensus}. For engineers and physicists working with such applications it is important to know that the global behaviour of the Kuramoto model on the Stiefel manifold is qualitatively different from that of the original Kuramoto model. For control applications, almost global synchronization is desirable since the probability of convergence does not decrease as $N$ increases. For model selection, the global behaviour of the real system should be taken into account.

\section{Problem Formulation}

\subsection{Notation}

The Frobenius inner product of $\ma{X}, \ma{Y}\in\R^{n\times p}$ is $g(\ma{X},\ma{Y})=\langle\ma{X},\ma{Y}\rangle=\trace\mat{X}\ma{Y}$. The norm of $\ma{X}$ is given by $\|\ma{X}\|=\smash{\langle\ma{X},\ma{X}\rangle^{\frac12}}$. The gradient on $\St\subset\R^{n\times p}$ (in terms of $g$) of a function $V:\St\rightarrow\R$ is given by $\nabla V=\Pi\widebar{\nabla} \widebar{V}$, where $\Pi:\R^{n\times p}\rightarrow\ts[\St]{\ma{X}}$ is an orthogonal projection operator, $\widebar{\nabla}$ denotes the gradient in the ambient Euclidean space,  and $\widebar{V}$ is any smooth extension of $V$ on $\R^{n\times p}$. 


A graph $\mathcal{G}$ is a pair $(\V,\E)$ where $\V=\{1,\ldots,N\}$ and $\E$ is a set of 2-element subsets of $\mathcal{V}$. Throughout this paper, if an expression depends on an edge $e\in\E$ and two nodes $i,j\in\V$, then it is implicitly understood that $e=e(i,j)=\{i,j\}$. Each element $i\in\V$ corresponds to a unique agent. Items associated with agent $i$ carry the subindex $i$; we let $\ma[i]{S}\in\St$ denote the state of an agent, $\Pi_i$ the orthogonal projection operator onto the tangent space $\ts[\St]{i}$ at $\ma[i]{S}$, $\mathcal{N}_i=\{j\in\V\,|\,\{i,j\}\in\E\}$ the neighbor set of $i$, $\nabla_i V$ the gradient of $V$ with respect to $\ma[i]{S}\in\St$, \etc

\subsection{The Stiefel manifold}

The compact, real Stiefel manifold $\St$ is the set of $p$-frames in $n$-dimensional Euclidean space $\R^n$ \citep{edelman1998geometry}. It can be embedded in $\R^{n\times p}$ as an analytic matrix manifold given by
\begin{align*}
\St=\{\ma{S}\in\R^{n\times p}\,|\,\mat{S}\ma{S}=\ma[p]{I}\}.
\end{align*}
The dimension of $\St$ is $np-\tfrac12p(p+1)$ due to the constraints. Important instances of Stiefel manifolds include the $n$-sphere $\mathsf{S}^n=\mathsf{St}(1,n+1)$, the special orthogonal group $\SO\simeq\mathsf{St}(n-1,n)$, and the orthogonal group $\mathsf{O}(n)=\mathsf{St}(n,n)$. Since $\|\ma{S}\|^2=p$ for all $\ma{S}\in\St$, it holds that $\St$ is a subset of the sphere of radius $\smash{p^\frac12}$ in the space of real $n\times p$ matrices. As rough guideline, the Stiefel manifold can be used to model systems whose states are constant in norm and subject to orthogonality constraints. 



Define the projections $\skews:\R^{n\times n}\rightarrow\so:\ma{X}\mapsto\tfrac{1}{2}(\ma{X}-\mat{X})$ and $\syms:\R^{n\times n}\rightarrow\so^\perp:\ma{X}\mapsto\tfrac12(\ma{X}+\mat{X})$. The tangent space of $\St$ at $\ma{S}$ is given by
\begin{align*}
\ts[\St]{\ma{S}}&=\{\ma{\Delta}\in\R^{n\times p}\,|\,\syms\mat{S}\!\ma{\Delta}=\ma{0}\}.
\end{align*}
Denote the tangent bundle of $\St$ by 
\begin{align*} 
\ts[\St]{}=\{(\ma{S},\ma{\Delta})\in\St\times \ts[\St]{\ma{S}}\}.
\end{align*}
The projection onto the tangent space, $\Pi:\St\times\R^{n\times p}\rightarrow\ts[\St]{\ma{S}}$, is given by 
\begin{align*}
\Pi(\ma{S},\ma{X})=\ma{S}\skews\mat{S}\!\ma{X}+(\ma[n]{I}-\ma{S}\mat{S})\ma{X}.
\end{align*}

\subsection{Synchronization on the Stiefel manifold}

The synchronization set, or consensus manifold, $\mathcal{C}$ of the $N$-fold product of a Stiefel manifold is defined as
\begin{align}\label{eq:C}
\mathcal{C}&=\{(\ma[i]{S})_{i=1}^N\in\St^N\,|\,\ma[i]{S}=\ma[j]{S},\forall\,\{i,j\}\in\E\},
\end{align}
where $(\ma[i]{S})_{i=1}^N$ denotes an $N$-tuple. The synchronization set is a (sub)manifold; it is diffeomorphic to $\St$ by the map $(\ma[i]{S})_{i=1}^N\mapsto \ma[1]{S}$. Let $d_{ij}=\|\ma[i]{S}-\ve[j]{S}\|$ be the chordal distance between agent $i$ and $j$. Given a graph $(\V,\E)$, define the potential function $V:\St^N\rightarrow\R$ by
\begin{align}
V&=\sum_{e\in\E}a_{ij}d_{ij}^2=\sum_{e\in\E}a_{ij}\|\ma[i]{S}-\ma[j]{S}\|^2\nonumber\\
&=2\sum_{e\in\E}a_{ij}(p-\langle\ma[i]{S},\ma[j]{S}\rangle),\label{eq:VSt}
\end{align}
where $a_{ij}\in(0,\infty)$ satisfies $a_{ij}=a_{ji}$ for all $e\in\E$. Note that $V$ is a real-analytic function, $V\geq0$, and $V|_{\mathcal{C}}=0$.

Denote $\ma{S}=(\ma[i]{S})_{i=1}^N$. Let $\widebar{V}:(\R^{n\times p})^N\rightarrow[0,\infty)$ be a  smooth extension of $V$ obtained by relaxing the requirement $\ma{S}\in\St^N$ to $\ma{S}\in(\R^{n\times p})^N$. We only need $\widebar{V}$ to define the gradient of $V$ in the embedding space $(\R^{n\times p})^N$ when restricted to $\St^N$. All smooth extensions hence give the same gradient \citep{loring2010an}. The system we study is the gradient descent flow on $\St^N$ given by
\begin{align}
\md{S}&=(\md[i]{S})_{i=1}^N=-\nabla V=(-\nabla_i V)_{i=1}^N,\nonumber\\
\md[i]{S}&=-\nabla_iV=-\Pi_i\widebar{\nabla}_i \widebar{V}=\Pi_i\sum_{j\in\Ni}a_{ij}\ma[j]{S}\label{eq:stateeq}\\
&=\ma[i]{S}\skews\Bigl(\mat[i]{S}\!\sum_{j\in\Ni}\!a_{ij}\ma[j]{S}\Bigr)+(\ma[n]{I}-\ma[i]{S}\mat[i]{S})\!\sum_{j\in\Ni}\!a_{ij}\ma[j]{S},\nonumber
\end{align}
where $\ma[i]{S}(0)\in\St$. Note that any equilibrium of \eqref{eq:stateeq} is a critical point of $V$ and vice versa.

Since the system \eqref{eq:stateeq} is an analytic gradient descent, it will converge to an equilibrium point from any initial condition \citep{lageman2007convergence}. This property allows us to adopt a strong definition of what it means for \eqref{eq:stateeq} to reach consensus:
\begin{definition}
The agents are said to synchronize, or to reach consensus, if $\lim_{t\rightarrow\infty}\ma{S}(t)\in\mathcal{C}$, where $\ma{S}$ is the state variable of the gradient descent flow \eqref{eq:stateeq} and $\mathcal{C}$ is the consensus manifold defined by \eqref{eq:C}.
\end{definition}

\subsection{Problem statement}



The aim of this paper is classify each instance of $\St$ as satisfying or not satisfying the following requirement: the gradient descent flow \eqref{eq:stateeq} with interaction topology given by any connected graph converges to the consensus manifold $\mathcal{C}$ from almost all initial conditions.


\subsection{High-dimensional Kuramoto model}

\label{sec:highdim}

We chose to define the high-dimensional Kuramoto model in complex networks over the Stiefel manifold $\St$ as
\begin{align}
\md[i]{X}=\ma[i]{\Omega}\ma[i]{X}+\ma[i]{X}\ma[i]{\Xi}-\nabla_i V, \quad\forall\,i\in\V,\label{eq:homo}
\end{align}
where $\ma[i]{X}\in\St$, $\ma[i]{\Omega}\in\so$, and $\ma[i]{\Xi}\in\mathsf{so}(p)$. The definition of \eqref{eq:homo} is motivated by two reasons as we detail in the next paragraphs. Note that \eqref{eq:homo} is a first-order model where the right-hand side is the sum of a drift-term and a gradient descent flow, just like for the Kuramoto model. The variables $\ma[i]{\Omega}$ and  $\ma[i]{\Xi}$ are  generalizations of the frequency term in the Kuramoto model. The expression  $\ma[i]{\Omega}\ma[i]{X}+\ma[i]{X}\ma[i]{\Xi}$ is not the standard form of an element of $\ts[\St]{i}$, but varying $\ma[i]{\Omega}$ and $\ma[i]{\Xi}$ spans the tangent space at any given $\ma[i]{X}$.

The model \eqref{eq:homo} encompasses the Kuramoto model. Better still, the following models are special cases of \eqref{eq:homo}:
\begin{align}
\md[i]{R}&=\ma[i]{\Omega}\ma[i]{R}+\!\sum_{j\in\Ni}\!a_{ij}\ma[i]{R}\skews\mat[i]{R}\ma[j]{R},\, \ma[i]{R}\in\SO,\label{eq:SO}\\
\vd[i]{x}&=\ma[i]{\Omega}\ve[i]{x}+(\ma[n+1]{I}-\ve[i]{x}\vet[j]{x})\!\sum_{j\in\Ni}\!a_{ij}\ve[j]{x},\, \ve[i]{x}\in\mathsf{S}^n,\label{eq:sphere}\\
\dot{\vartheta}_i&=\omega_i+\sum_{j\in\Ni}a_{ij}\sin(\vartheta_j-\vartheta_i),\,\vartheta_i\in\R,\label{eq:circle}
\end{align}
where $\ma[i]{\Omega}\in\so$, and $\omega_i\in\R$, and each system consists of $N$ equations; one for each $i\in\V$. 


To get \eqref{eq:SO} from \eqref{eq:homo} , let $p=n$ and set $\ma[i]{R}=\ma[i]{X}$,  $\ma[i]{\Xi}=\ma{0}$. Note that $\Pi_i:\R^{n\times n}\rightarrow\ts[\mathsf{O}(n)]{i}$ is given by $\Pi_i\ma{Y}=\ma[i]{X}\skews\mat[i]{X}\ma{Y}$ since $\ma[i]{X}\mat[i]{X}=\ma[n]{I}$. The restriction of $\ma[i]{R}(0)\in\SO$ implies that $\ma[i]{R}(t)\in\SO$ for all $t\in[0,\infty)$. To get \eqref{eq:sphere} from \eqref{eq:homo}, let $p=1$ and set $\ve[i]{x}=\ma[i]{X}$. Note that $\Pi_i:\R^{n+1\times 1}\rightarrow\ts[\mathsf{S}^n]{i}$ is given by $\Pi_i\ma{y}=(\ma[n+1]{I}-\ve[i]{x}\vet[i]{x})\ve[i]{y}$. To get \eqref{eq:circle} from \eqref{eq:sphere} (and hence also from \eqref{eq:homo} via \eqref{eq:sphere}), let $n=2$, $\ve[i]{x}=[\cos\vartheta_i\,\sin\vartheta_i]\mtr$,  $\omega_i=\langle\ve[2]{e},\ma[i]{\Omega}\ve[1]{e}\rangle$ and solve for $\dot{\vartheta}_i$.


The cases of homogeneous frequencies and zero frequencies are equivalent; \ie  \eqref{eq:homo} is equivalent to \eqref{eq:stateeq} in the case of $\ma[i]{\Omega}=\ma{\Omega}$, $\ma[i]{\Xi}=\ma{\Xi}$. To see this, introduce the variables $\ma{R}=\exp(-t\ma{\Omega})\in\SO$, $\ma{Q}=\exp(-t\ma{\Xi})\in\mathsf{SO}(p)$, form a rotating coordinate frame $\ve[i]{S}=\ma{R}\ma[i]{X}\ma{Q}\in\St$, and change variables
\begin{align*}
\md[i]{S}={}&-\ma{R}\ma{\Omega}\ma[i]{X}\ma{Q}+\ma{R}\md[i]{X}\ma{Q}-\ma{R}\ma[i]{X}\ma{\Xi}\ma{Q}\\
={}&-\ma{R}\nabla_i V(\ma[i]{X})_{i=1}^N\ma{Q}\\
={}&\ma{R}\ma[i]{X}\ma{Q}\mat{Q}\skews\Bigl(\mat[i]{X}\mat{R}\ma{R}\sum_{j\in\Ni}a_{ij}\ma[j]{X}\Bigr)\ma{Q}+\\
&\ma{R}(\ma[n]{I}-\ma[i]{X}\ma{Q}\mat{Q}\mat[i]{X})\mat{R}\ma{R}\sum_{j\in\Ni}a_{ij}\ma[j]{X}\ma{Q}\\
={}&\ma[i]{S}\skews\Bigl(\mat[i]{S}\sum_{j\in\Ni}a_{ij}\ma[j]{S}\Bigr)+(\ma[n]{I}-\ma[i]{S}\mat[i]{S})\sum_{j\in\Ni}a_{ij}\ma[j]{S}.
\end{align*}

\subsection{Local stability and global attractiveness}

The results of this paper concern the global stability properties of the flow \eqref{eq:stateeq}. The local stability properties of the system are summarized in Proposition \ref{th:stability}. This result states some rather generic properties of analytic gradient descent flows. We do not give a proof, but refer the interested reader to \cite{lageman2007convergence,helmke2012optimization}.

\begin{proposition}\label{th:stability}
The gradient descent flow \eqref{eq:stateeq} converges to a critical point of $V$. The sublevel sets 
	\begin{align*}
	\mathcal{L}(h)=\{\ma{S}\in\St^N\,|\,V(\ma{S})\leq h\}
	\end{align*}
are forward invariant.
\end{proposition}

Note that all global minimizers of $V$ belong to $\mathcal{C}$ since $V\geq0$ with equality only if $\ma{S}\in\mathcal{C}$. From $\dot{V}=\langle\nabla V,\md{S}\rangle=-\|\nabla V\|^2$ it follows that $\mathcal{C}$ is stable. Let $\mathcal{Q}$ denote all critical points of $V$ that are disjoint from $\mathcal{C}$. The distance between $\mathcal{C}$ and $\mathcal{Q}$ is positive, wherefore $\mathcal{C}$ is asymptotically stable. By Proposition \ref{th:stability}, the region of attraction of $\mathcal{C}$ contains the largest sublevel set $\mathcal{L}(h)$ which is disjoint from $\mathcal{Q}$.

\begin{definition}
An equilibrium set $\mathcal{Q}\subset\St^N$ of system \eqref{eq:stateeq} is referred to as almost globally asymptotically stable (\AGAS) if it is stable and attractive from all initial conditions $\ma{S}(0)\in\smash{\St^N}\backslash\mathcal{N}$, where $\mathcal{N}\subset\smash{\St^N}$ has Haar measure zero on $\smash{\St^N}$. 
\end{definition}

It is not possible to globally stabilize an equilibrium set on a compact manifold by means of continuous, time-invariant feedback \citep{bhat2000topological}. This obstruction, which is due to topological reasons, does not exclude the possibility of a set being \AGAS{}.

\section{Main Result}
\label{sec:main}


\begin{theorem}\label{th:main}
Let the pair $(p,n)$ satisfy $p\leq \tfrac{2}{3}n-1$ and $\mathcal{G}$ be connected. The consensus manifold 
\begin{align*}
\mathcal{C}&=\{(\ma[i]{S})_{i=1}^N\in\St^N\,|\,\ma[i]{S}=\ma[j]{S},\forall\,\{i,j\}\in\E\},
\end{align*}
is an \AGAS{} equilibrium set of the gradient descent flow on $\St^N$ given by
\begin{align*}
\md[i]{S}&=\ma[i]{S}\skews\Bigl(\mat[i]{S}\!\sum_{j\in\Ni}\!a_{ij}\ma[j]{S}\Bigr)+(\ma[n]{I}-\ma[i]{S}\mat[i]{S})\!\sum_{j\in\Ni}\!a_{ij}\ma[j]{S}.
\end{align*}
%
\end{theorem}

The calculations involved in the proof of Theorem \ref{th:main} are extensive. We give a brief proof sketch that covers the main ideas. All the details are provided in Appendix \ref{sec:critical} to \ref{sec:nlp}. 

\begin{pf}
If the linearization of \eqref{eq:stateeq} around an equilibrium $\ma{S}=(\ma[i]{S})_{i=1}^N\in\St$ has an eigenvalue with strictly positive real part, then that equilibrium is exponentially unstable by the indirect method of Lyapunov. We can also think of equilibria as critical points of $V$, \ie points where the gradient is the zero vector. The nature of a critical point can often be determined by studying the Riemannian Hessian $\ma{H}(\ma{S})$ of $V$, \ie the first non-zero term in the Taylor expansion of $V$. Note that the Hessian matrix equals the linearization matrix, albeit multiplied by minus one. The instability criterion given by the indirect method of Lyapunov is hence equivalent to the necessary second-order optimality conditions.
	
Any set of exponentially unstable equilibria of a pointwise convergent system have a measure zero region of attraction \citep{freeman2013global}. Pointwise convergence means, roughly speaking, that the system does not admit any limit cycles. Every trajectory converges to some point. Gradient descent flows of analytic functions on compact analytic manifolds are pointswise convergent as a consequence of the \L{}ojasiewicz gradient inequality \citep{lageman2007convergence}. The consensus manifold $\mathcal{C}$ is stable by Lyapunov's theorem since 
$\dot{V}=\langle\nabla V,\vd{S}\rangle=-\|\nabla V\|^2$. It follows that $\mathcal{C}$ is $\AGAS$ if $\ma{H}(\ma{S})$ evaluated at any equilibrium $\ma{S}\notin\mathcal{C}$ has an eigenvalue with strictly negative real part. 

Let $q:\ts[\St^N]{}\rightarrow\R$ denote the quadratic form obtained from the Riemannian Hessian $\ma{H}(\ma{S})$ evaluated at a critical point $\ma{S}\in\St^N$. The Hessian at $\ma{S}\in\St^N$ is a symmetric linear operator $\ma{H}:\ts[\St^N]{\ma{S}}\rightarrow\ts[\St^N]{\ma{S}}$  in the sense that
\begin{align*} \langle(\ma[i]{X})_{i=1}^N,\ma{H}(\ma{S})(\ma[i]{Y})_{i=1}^N\rangle=\langle\ma{H}(\ma{S})(\ma[i]{X})_{i=1}^N,(\ma[i]{Y})_{i=1}^N\rangle
\end{align*}
\citep{absil2009optimization}. As such, its eigenvalues are real. The quadratic form $q$ therefore bounds the smallest  eigenvalue of the linear operator $\ma{H}(\ma{S})$ from above.
%
%
Our goal is to establish exponential instability of  all equilibria $\ma{S}\notin\mathcal{C}$ by finding a tangent vector  $(\ma[i]{\Delta})_{i=1}^N\in\ts[\St^N]{\ma{S}}$ such that 
\begin{align*}
q((\ma[i]{S})_{i=1}^N,(\ma[i]{\Delta})_{i=1}^N)=\langle(\ma[i]{\Delta})_{i=1}^N,\ma{H}(\ma{S})(\ma[i]{\Delta})_{i=1}^N\rangle<0.
\end{align*}
We want to use a tangent vector $(\ma[i]{\Delta})_{i=1}^N$ whose representation in the eigenvector basis of $\ts[\St^N]{\ma{S}}$ is dominated by the eigenvector of $\ma{H}(\ma{S})$ with the smallest eigenvalue. The quadratic form $q$ will then approximate the smallest eigenvalue multiplied by $\|(\ma[i]{\Delta})_{i=1}^N\|^2$.

Consider tangent vectors pointing towards $\mathcal{C}$, \ie $\ma[i]{\Delta}=\Pi_i\ma{\Delta}$ for some $\ma{\Delta}\in\R^{n\times p}$. The intuition for this choice is that a small perturbation of the system where every agent is moved in the same direction should not result in an increase of $V$ (if the perturbations are similar they cancel each other for each pair $(i,j)\in\St$). Moreover, it is possible that there is a net increase in cohesion which would yield a decrease in $V$. We do not need to find an expression for the desired tangent vector, it suffices to prove that it exists. 

We show that $q$ only assumes negative values by solving an optimization problem to minimize an upper bound of $q$ over $\ts[\St^N]{}$. The upper bound is obtained by relaxing the complex network of relations between agents at an equilibrium and only consider the effect of pairwise interactions. For any equilibrium $\ma{S}\notin\mathcal{C}$ and pair $(p,n)$ such that $p\leq\tfrac{2n}{3}-1$, we find that there is a tangent vector towards $\mathcal{C}$ which results in the upper bound on $q$ being strictly negative. Any equilibrium $\ma{S}\notin\mathcal{C}$ is hence exponentially unstable. Throughout these steps, we do not utilize any particular property of the graph topology except connectedness.\hfill\phantom{.}\qed\end{pf}

\begin{remark}
The inequality $p\leq \tfrac23n-1$ is sufficient for $\mathcal{C}$ to be \AGAS. In a more general setting of Kuramoto models on closed Riemannian manifolds, it can be showed that a manifold being multiply connected precludes $\mathcal{C}$ being \AGAS{}. A multiply connected manifold is, roughly speaking, a manifold with a hole, for example a torus. In particular, the only multiply connected Stiefel manifolds are $\mathsf{St}(n-1,n)\simeq\SO$ and $\mathsf{St}(n,n)=\mathsf{O}(n)$  \citep{james1976topology}. Further results on multistability of the Kuramoto model on $\SO=\{\ma{S}\in\mathsf{St}(n,n)\,|\,\det\ma{S}=1\}$ are given in \citet{deville2018synchronization}. The question if $\mathcal{C}$ is \AGAS{} for all connected graphs on $\St$ where $\tfrac23n-1<p\leq n-2$ remains open. Using Monto Carlo experiments to estimate the probability measure of the region of attraction of $\mathcal{C}$, we observe that $\mathcal{C}$ appears to be \AGAS{} on some such Stiefel manifolds for networks over which $\mathsf{St}(n-1,n)$ is multistable.
\end{remark}

\section{Numerical Examples}
\label{sec:numerical}

\begin{figure}[htb!]
	\centering	\includegraphics[trim={1.3cm 1.3cm 1.3cm 1.3cm}, clip, width=0.2\textwidth]{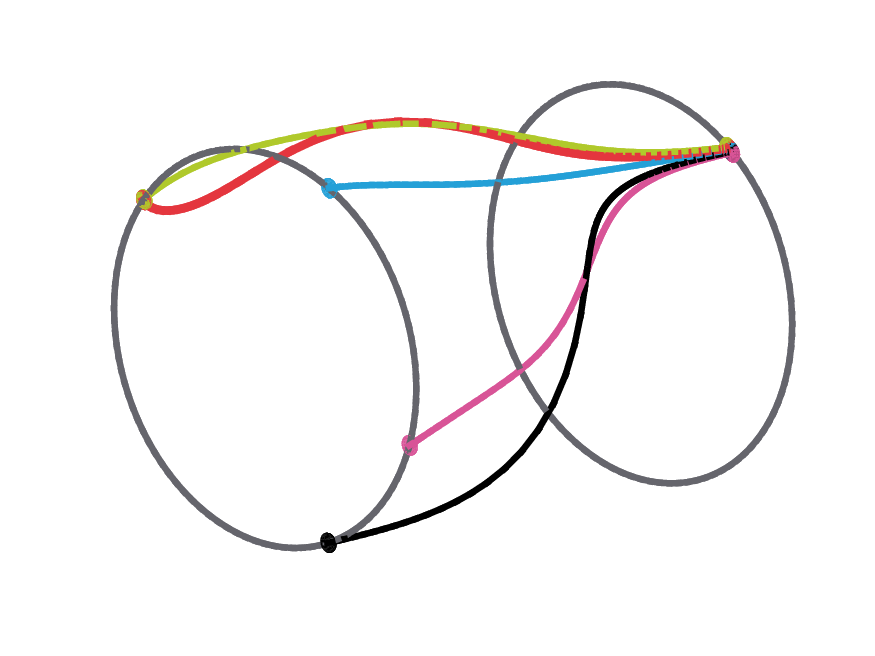}\hspace{1cm}
	\includegraphics[trim={1.3cm 1.3cm 1.3cm 1.3cm}, clip, width=0.2\textwidth]{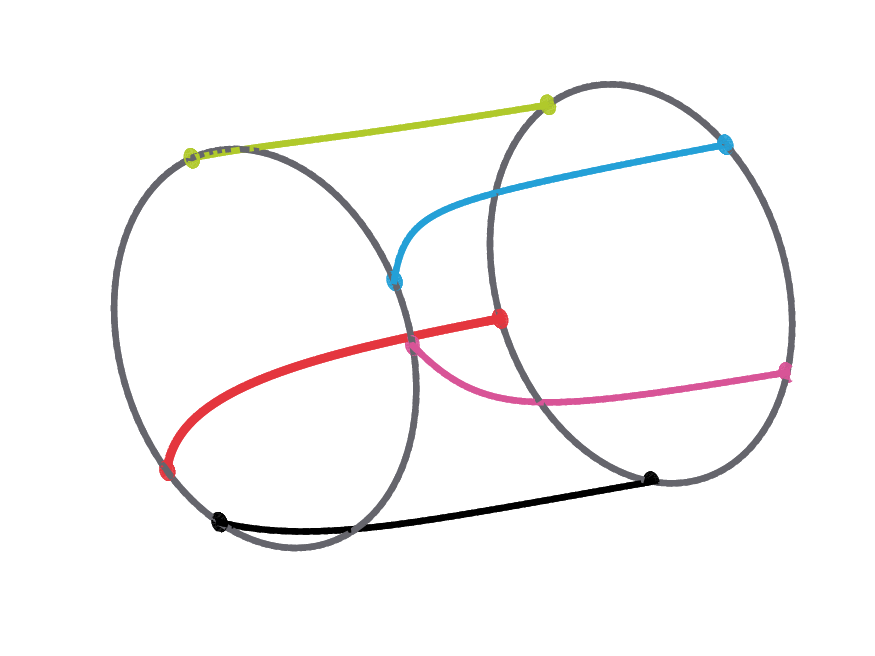}
	\caption{Two sets of trajectories for five agents on $\GC$ that are connected by the graph $\mathcal{H}_5$. The agents evolve from random initial conditions towards the sets $\mathcal{C}$ (left) and $\mathcal{Q}_{12}$ (right). The positive direction of time is from left to right in both figures.}
	\label{fig:circle}
\end{figure}

\begin{figure}[htb!]
	\centering	
	\includegraphics[trim={1.5cm 1.5cm 1.5cm 1.5cm}, clip, width=0.3\textwidth]{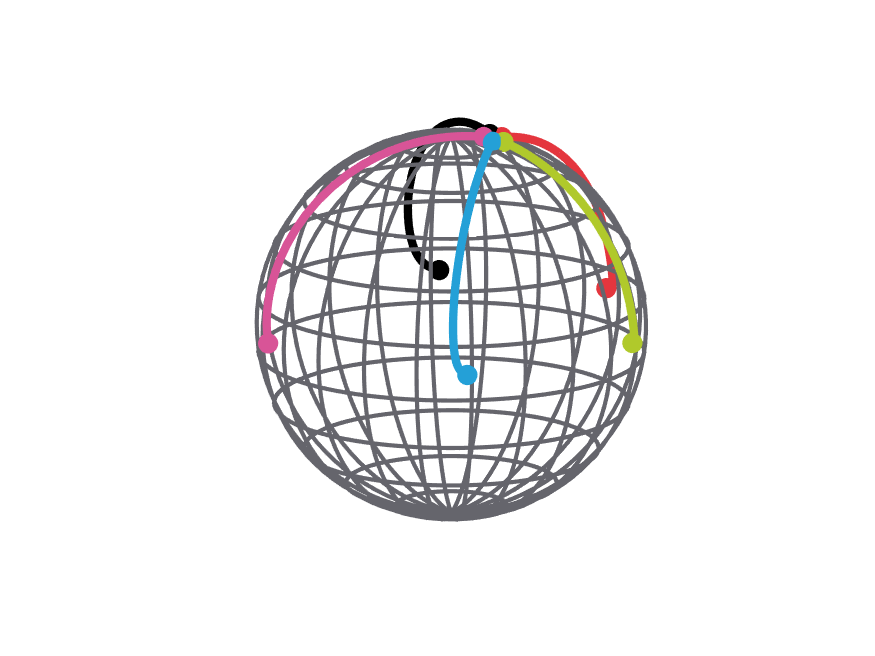}
	\caption{The trajectories of five agents with on $\twosphere$ that are connected by the graph $\mathcal{H}_5$. The agents evolve from a point  close to $\mathcal{Q}_{13}$ (\ie close to the equator) towards $\mathcal{C}$ near the north pole.}
	\label{fig:sphere}
\end{figure}

We provide numerical examples to illustrate the evolution of system \eqref{eq:stateeq} on $\mathsf{St}(1,2)=\mathcal{S}^1$, $\mathsf{St}(1,3)=\mathcal{S}^2$, and $\mathsf{St}(2,3)\simeq\SOT$ when $a_{ij}=1$. Let $\mathcal{H}_N$ denote the cyclic graph over $N$ nodes, \ie
\begin{align*}
\mathcal{H}_N=(\{1,\ldots,N\},\{\{i,j\}\subset\V\,|\,j=i+1 \}),
\end{align*}
where we set $N+1=1$. The equilibrium set
\begin{align*}
\mathcal{Q}_{1n}=\{&(\ma[i]{x})_{i=1}^N\in(\Sn)^N\,|\,\exists\, \ma{R}\in\mathsf{SO}(n),\\
&\tfrac{1}{\sqrt2}\|\Log\ma{R}\|=\tfrac{2\pi}{N},\,\ma[i+1]{x}=\ma{R}\ma[i]{x},\,\forall\,i\in\V\},
\end{align*}
is asymptotically stable for the system \eqref{eq:stateeq} if $n=1$ and  $N\geq5$, but unstable for all $N\in\N$ if $n\geq2$. This is illustrated in Fig. \ref{fig:circle} and \ref{fig:sphere}.

To understand this difference, note that the complement of the circle is two open hemispheres. The consensus manifold $\mathcal{C}$ is asymptotically stable on any open hemisphere \citep{markdahl2018tac}. As such, we may move each agent an arbitrarily small distance from $\mathcal{Q}_{13}$, perturbing them into an open hemisphere, whereby they will reach consensus. 

\begin{figure}[htb!]
	\centering	\includegraphics[trim={0.8cm 0.8cm 0.8cm 0.8cm}, clip, width=0.3\textwidth]{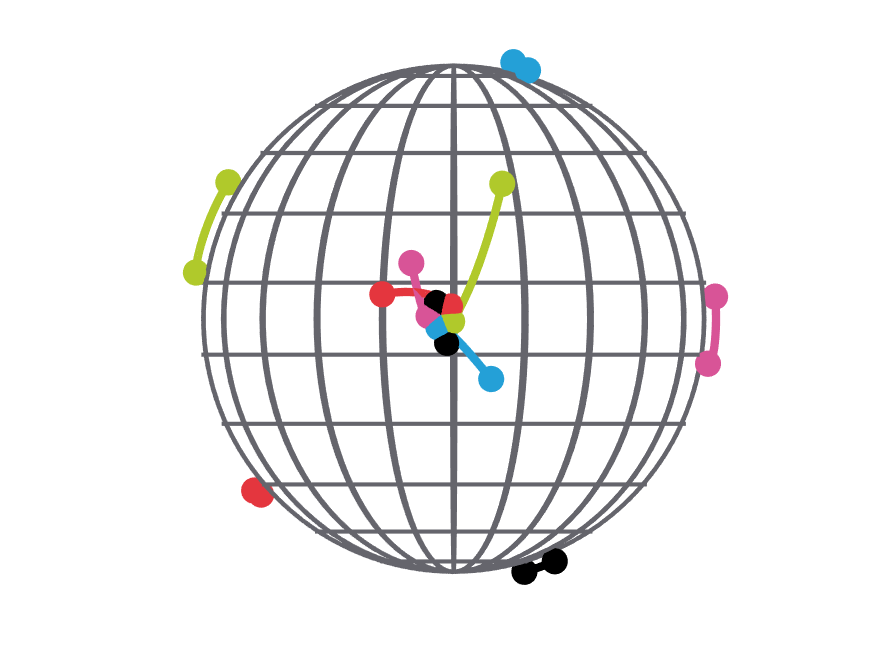}
	\caption{The trajectories of five agents on $\mathsf{St}(2,3)$ that are connected by the graph $\mathcal{H}_5$. Each agent state is represented as an orthogonal pair of vectors on $\twosphere$. The agents are initially perturbed away from the equilibrium set $\mathcal{Q}_{23}$ but ultimately end up close to it.}
	\label{fig:st23}
\end{figure}

Each element of $\mathsf{St}(2,3)$ is a pair of orthogonal unit vectors $(\ma[i]{S}\ve[1]{e},\ma[i]{S}\ve[2]{e})\in\mathcal{S}^2\times\mathcal{S}^2$. They can be visualized as pairs of points on a single sphere. Consider the equilibrium set
\begin{align*}
\mathcal{Q}_{23}=\{&(\ma[i]{S})_{i=1}^N\in\StN|\,\ma[i+1]{S}\ve[1]{e}\!=\ma[i]{S}\ve[1]{e},\\
&\exists\,\ma{R}\in\SOT,\,\tfrac{1}{\sqrt2}\|\Log\ma{R}\|=\tfrac{2\pi}{N},\\
&\ma[i+1]{S}=\ma{R}\ma[i]{S},\,\forall\,i\in\V\}
\end{align*}
on $\mathsf{St}(2,3)\simeq\SOT$. In $\mathcal{Q}_{23}$, the first unit vectors $\ma[i]{S}\ve[1]{e}$ are aligned with each other while the second unit vectors $\ma[i]{S}\ve[2]{e}$ are spread out over a great circle. If the states are slightly perturbed to leave $\mathcal{Q}_{23}$, then they will often stay close to $\mathcal{Q}_{23}$ for all future times, see Fig. \ref{fig:st23}.

Note the difference in behavior of system \eqref{eq:stateeq} on $\mathcal{S}^2$ and $\SOT$. Why does the high-dimensional system on $\mathcal{S}^2$ reach consensus while the system on $\mathsf{St}(p,n)$ does not? Roughly speaking, the first vectors $\ma[i]{S}\ve[1]{e}$ all remain close to each other and this constrains the second  vectors $\ma[i]{S}\ve[2]{e}$ to a tubular neighborhood of the great circle they started out on. The dynamics on the tubular neighborhood are sufficiently similar to the Kuramoto model on the circle that the second unit vectors ultimately converge to a configuration that is similar to $\mathcal{Q}_{12}$ in Fig. \ref{fig:circle}.

%

\section{Conclusions and Future Work}

This paper formulates a Kuramoto model on the Stiefel manifold and studies its global behaviour. The Stiefel manifold includes both instances on which synchronization is multistable, \ie the Kuramoto model on the circle and the Lohe model on the special orthogonal group $\SO$ \citep{deville2018synchronization}, and instances on which synchronization is almost globally stable, \ie the $n$-sphere for $n\in\N\backslash\{1\}$ \citep{markdahl2018tac}. As such, studying its global behaviour can give us further insight into the global behaviour of consensus seeking systems on more general manifolds. The consensus manifold on $\St$ is \AGAS{} if the pair $(p,n)$ satisfies $p\leq\tfrac23n-1$. We believe that this condition is conservative due to the inequalities involved in calculating an upper bound on the smallest eigenvalue of the Riemannian Hessian, see Appendix \ref{sec:upper} and \ref{sec:nlp}. Rather, we conjecture that a sharp inequality is given by $p\leq n-2$, corresponding to all the simply connected Stiefel manifolds \citep{james1976topology}. Related topics will be explored in future work. 

\section{Acknowledgments}

The authors would like to thank the anonymous reviewers.

\bibliographystyle{named}        
\bibliography{autosam}           

\appendix

\section{Appendix}
\label{app:app}

\subsection{Equilibria are critical points}
\label{sec:critical}

We start by characterizing the equilibria of system \eqref{eq:stateeq}. At an equilibrium,
\begin{align*}
\ma[i]{S}\skews\Bigl(\mat[i]{S}\sum_{j\in\Ni}a_{ij}\ma[j]{S}\Bigr)+(\ma[n]{I}-\ma[i]{S}\mat[i]{S})\sum_{j\in\Ni}a_{ij}\ma[j]{S}=\ma{0}.
\end{align*}
Since the two terms in this expression are orthogonal, we get
\begin{align}
\begin{split}
\skews\Bigl(\mat[i]{S}\sum_{j\in\Ni}a_{ij}\ma[j]{S}\Bigr)&=\ma{0},\\
(\ma[n]{I}-\ma[i]{S}\mat[i]{S})\sum_{j\in\Ni}a_{ij}\ma[j]{S}&=\ma{0}.\label{eq:eq}
\end{split}
\end{align}
Assume \eqref{eq:eq} holds. Define $\ve[i]{\Sigma}=\sum_{j\in\Ni}a_{ij}\ma[j]{S}$. Since $\ma[i]{\Sigma}=\ma[i]{S}\mat[i]{S}\ma[i]{\Sigma}$, it follows that $\ma[i]{\Sigma}\in\im\ma[i]{S}$. Hence $\ve[i]{\Sigma}=\ma[i]{S}\ma[i]{\Gamma}$ for some $\ma[i]{\Gamma}\in\R^{p\times p}$. Moreover, since $\skews\mat[i]{S}\ma[i]{\Sigma}=\skews\ma[i]{\Gamma}=\ma{0}$, we find that $\ma[i]{\Gamma}$ is symmetric. 

\subsection{The Hessian on $\St^N$}

The next step in the proof sketch of Theorem \ref{th:main} is to determine the Hessian  $\ma{H}=[\nabla_k(\nabla_i V)_{st}]$. Let $\maw[i,st]{F}=(\Pi_i\widebar{\nabla}_i\widebar{V})_{st}:\R^{N\times n\times p}\rightarrow\R$ be a smooth extension of $\ma[i,st]{F}=(\nabla_i V)_{st}=\langle\ve[s]{e},\nabla_i V\ve[t]{e}\rangle:\St^N\rightarrow\R$ obtained by relaxing the constraint $\ma[i]{S}\in\St$ to $\ma[i]{S}\in\R^{n\times p}$. Take a $k\in\V$ and calculate
\begin{align*}
\widebar{\nabla}_k\maw[i,st]{F}={}&\widebar{\nabla}_k(\Pi_i\widebar{\nabla}_i\widebar{V})_{st}=\widebar{\nabla}_k\langle\ve[s]{e},\Pi_i\widebar{\nabla}_i\widebar{V}\ve[t]{e}\rangle\\
={}&\widebar{\nabla}_k\Bigl\langle\ve[s]{e},\Bigl(-\ma[i]{S}\skews\Bigl(\mat[i]{S}\sum_{j\in\Ni}a_{ij}\ma[j]{S}\Bigr)-\\
&\hspace{1.4cm}(\ma[n]{I}-\ma[i]{S}\mat[i]{S})\sum_{j\in\Ni}a_{ij}\ma[j]{S}\Bigr)\ve[t]{e}\Bigr\rangle\\
={}&-\widebar{\nabla}_k\Bigl\langle\ve[s]{e},\ma[i]{S}\skews\Bigl(\mat[i]{S}\sum_{j\in\Ni}a_{ij}\ma[j]{S}\Bigr)\ve[t]{e}\Bigr\rangle-\\
&\widebar{\nabla}_k\Bigl\langle\ve[s]{e},\sum_{j\in\Ni}a_{ij}\ma[j]{S}\ve[t]{e}\Bigr\rangle+\\
&\widebar{\nabla}_k\Bigl\langle\ve[s]{e},\ma[i]{S}\mat[i]{S}\sum_{j\in\Ni}a_{ij}\ma[j]{S}\ve[t]{e}\Bigr\rangle.
\end{align*}	
Using the rules governing derivatives of inner products with respect to matrices, introducing  $\ma[st]{E}=\ve[s]{e}\vet[t]{e}=\ve[s]{e}\otimes\ve[t]{e}$, after a few calculations, we obtain
\begin{align*}	
\widebar{\nabla}_k\maw[i,st]{F}={}&\begin{cases}
-a_{ik}\Pi_i\ma[st]{E} & \textrm{ if } k\in\Ni,\\
\ma[st]{E}\skews\left(\mat[i]{S}\sum_{j\in\Ni}a_{ij}\ma[j]{S}\right)+&\\
\sum_{j\in\Ni}a_{ij}\ma[j]{S}\syms(\mat[i]{S}\ma[st]{E})+& \\	
\ma[st]{E}\sum_{j\in\Ni}a_{ij}\mat[j]{S}\ma[i]{S}&\textrm{ if } k=i,\\
\ma{0}& \textrm{ otherwise.}
\end{cases}
\end{align*}
Evaluate at an equilibrium, where $\sum_{j\in\Ni}a_{ij}\ma[j]{S}=\ma[i]{S}\ma[i]{\Gamma}$ and $\ma[i]{\Gamma}\in\R^{p\times p}$ is symmetric by Section \ref{sec:critical}, to find
\begin{align*}
\widebar{\nabla}_k\maw[i,st]{F}={}&\begin{cases}
-a_{ik}\Pi_i\ma[st]{E} & \textrm{ if } k\in\Ni,\\
\ma[i]{S}\ma[i]{\Gamma}\syms(\mat[i]{S}\ma[st]{E})+\ma[st]{E}\ma[i]{\Gamma} & \textrm{ if } k=i,\\
\ma{0}& \textrm{ otherwise.}
\end{cases}
\end{align*}

The Hessian on $\St^N$ is a $(N\times n\times p)^2$-tensor consisting of $N^2np$ blocks $\ma[ki,st]{H}\in\R^{n\times p}$ formed by projecting the Hesssian in $(\R^{n\times p})^N$ on the tangent space of $\ma[k]{S}$
\begin{align*}
\ma[ki,st]{H}&=\nabla_{k}(\nabla_{i}V)_{st}=\Pi_k\widebar{\nabla}_{k}\maw[i,st]{F}\\
&=\Pi_k\widebar{\nabla}_{k}(\Pi_i\widebar{\nabla}_i\widebar{V})_{st}.
\end{align*}

\subsection{The quadratic form}

The quadratic form $q:\ts[\St^N]{}\rightarrow\R$ determines the nature of a critical point $\ma{S}$ in the sense of the necessary second-order optimality conditions \citep{nocedal1999numerical}. Consider the quadratic form obtained from the Hessian $\ma{H}(\ma{S})$ evaluated at an equilibrium $\ma{S}$ together with a tangent vector $(\ma[i]{\Delta})_{i=1}^N\in\ts[\St^N]{\ma{S}}$, where $\ma[i]{\Delta}=\Pi_i\ma{\Delta}$ for some $\ma{\Delta}\in\R^{n\times p}$, \ie the tangent vector is pointing towards the consensus manifold $\mathcal{C}$,
\begin{align*}
q={}&\sum_{i=1}^N\sum_{k=1}^N\langle\ma[i]{\Delta},[\langle\ma[k]{\Delta},\nabla_k(\nabla_i V)_{st}\rangle]\rangle\\
={}&\sum_{i=1}^N\sum_{k=1}^N\langle\Pi_i\ma{\Delta},[\langle\Pi_k\ma{\Delta},\Pi_k\widebar{\nabla}_k\maw[i,st]{F}\rangle]\rangle.
\end{align*}

Note that 
$\langle\Pi_k\ma{X},\Pi_k\ma{Y}\rangle=\langle\Pi_k\ma{X},\ma{Y}\rangle$. The quadratic form is hence
\begin{align*}
q={}&\sum_{i=1}^N\sum_{k=1}^N\langle\Pi_i\ma{\Delta},[\langle\Pi_k\ma{\Delta},\widebar{\nabla}_k\maw[i,st]{F}\rangle]\rangle.
\end{align*}
Denote $\ma[ki,st]{P}=\langle\Pi_k\ma{\Delta},\widebar{\nabla}_k\maw[i,st]{F}\rangle$. Then
\begin{align*}
\ma[ki,st]{P}=\begin{cases}
\langle\Pi_k\ma{\Delta},-a_{ik}\Pi_i\ma[st]{E}\rangle\\
\langle\Pi_i\ma{\Delta},\ma[i]{S}\ma[i]{\Gamma}\syms(\mat[i]{S}\ma[st]{E})+\ma[st]{E}\ma[i]{\Gamma}\rangle\\
0
\end{cases}
\end{align*}
for the cases of $k\in\Ni$, $k=i$, and $k\notin\Ni\cup\{i\}$ respectively. Denote $\ma[ki]{P}=[\ma[ki,st]{P}]$ and calculate
\begin{align*}
\ma[ki]{P}&=\begin{cases}
-a_{ik}\Pi_i\Pi_k\ma{\Delta} & \textrm{ if } k\in\Ni,\\
\ma[i]{S}\syms\mat[i]{\Sigma}\Pi_i\ma{\Delta}+\Pi_i(\ma{\Delta})\ma[i]{\Gamma} & \textrm{ if } k=i,\\
\ma{0} & \textrm{ otherwise.}
\end{cases}
\end{align*}
To see this, consider each case separately. For $k\in\Ni$,
\begin{align*}
\ma[ki,st]{P}={}&\langle(\ma[n]{I}-\Pi_i+\Pi_i)\Pi_k\ma{\Delta},-a_{ik}\Pi_i\ma[st]{E}\rangle\\
={}&-a_{ik}\langle\Pi_i\Pi_k\ma{\Delta},\ma[st]{E}\rangle=-a_{ik}(\Pi_i\Pi_k\ma{\Delta})_{st},
\end{align*}
whereby $\ma[ki]{P}=-a_{ik}\Pi_i\Pi_k\ma{\Delta}$. For the case of $k=i$,
\begin{align*}
\ma[ii,st]{P}={}&\langle\Pi_i\ma{\Delta},\ma[i]{S}\ma[i]{\Gamma}\syms(\mat[i]{S}\ma[st]{E})+\ma[st]{E}\ma[i]{\Gamma}\rangle\\
={}&\trace(\Pi_i\ma{\Delta})\mtr(\tfrac12\ma[i]{\Sigma}(\mat[i]{S}\ma[st]{E}+\mat[st]{E}\ma[i]{S})+\ma[st]{E}\ma[i]{\Gamma})\\
={}&\tfrac12\trace((\Pi_i\ma{\Delta})\mtr\ma[i]{\Sigma}\mat[i]{S}\ma[st]{E}+\ma[i]{S}(\Pi_i\ma{\Delta})\mtr\ma[i]{\Sigma}\mat[st]{E})+\\
&\trace\ma[i]{\Gamma}(\Pi_i\ma{\Delta})\mtr\ma[st]{E}\\
={}&\tfrac12(\ma[i]{S}\mat[i]{\Sigma}\Pi_i\ma{\Delta})_{st}+\tfrac12(\ma[i]{S}(\Pi_i\ma{\Delta})\mtr\ma[i]{\Sigma})_{st}+\\
&(\Pi_i(\ma{\Delta})\ma[i]{\Gamma})_{st},
\end{align*}
whereby $\ma[ii]{P}=\ma[i]{S}\syms\mat[i]{\Sigma}\Pi_i\ma{\Delta}+\Pi_i(\ma{\Delta})\ma[i]{\Gamma}$.

This gives us the quadratic form
\begin{align*}
q={}&\sum_{i=1}^N\sum_{k=1}^N\langle\Pi_i\ma{\Delta},[\ma[ki,st]{P}]\rangle=\sum_{i=1}^N\sum_{k=1}^N\langle\Pi_i\ma{\Delta},\ma[ki]{P}\rangle\\
={}&\sum_{e\in\mathcal{E}}\langle\Pi_i\ma{\Delta},\ma[ki]{P}\rangle+\langle\Pi_k\ma{\Delta},\ma[ki]{P}\rangle+\sum_{i\in\mathcal{V}}\langle\Pi_i\ma{\Delta},\ma[ii]{P}\rangle.
\end{align*}
For ease of notation, let $q=2\sum_{e\in\E}q_{ik}+\sum_{i\in\V}q_i$, where
\begin{align*}
q_{ik}={}&\langle\Pi_i\ma{\Delta},\ma[ki]{P}\rangle=-a_{ik}\langle\Pi_i\ma{\Delta},\Pi_k\ma{\Delta}\rangle=q_{ki},\\
q_i={}&\langle\Pi_i\ma{\Delta},\ma[ii]{P}\rangle.
\end{align*}
%


Calculate
\begin{align*}
q_{ik}={}&-a_{ik}\langle\Pi_i\ma{\Delta},\Pi_k\ma{\Delta}\rangle\\
={}&a_{ik}(-\langle\ma{\Delta},\ma{\Delta}\rangle+\tfrac12\langle\ma[i]{S}(\mat[i]{S}\ma{\Delta}+\mat{\Delta}\ma[i]{S}),\ma{\Delta}\rangle+\\
&\tfrac12\langle\ma{\Delta},\ma[k]{S}(\mat[k]{S}\ma{\Delta}+\mat{\Delta}\ma[k]{S})\rangle-\\
&\tfrac14\langle\ma[i]{S}(\mat[i]{S}\ma{\Delta}+\mat{\Delta}\ma[i]{S}),\ma[k]{S}(\mat[k]{S}\ma{\Delta}+\mat{\Delta}\ma[k]{S})\rangle)\\
={}&a_{ik}\trace(-\mat{\Delta}\ma{\Delta}+\tfrac12\mat{\Delta}\ma[i]{S}\mat[i]{S}\ma{\Delta}+\tfrac12\mat[i]{S}\ma{\Delta}\mat[i]{S}\ma{\Delta}+\\
&\hphantom{a_{ik}\trace(}\tfrac12\mat{\Delta}\ma[k]{S}\mat[k]{S}\ma{\Delta}+\tfrac12\mat{\Delta}\ma[k]{S}\mat{\Delta}\ma[k]{S}-\\
&\hphantom{a_{ik}\trace(}\tfrac14\mat{\Delta}\ma[i]{S}\mat[i]{S}\ma[k]{S}\mat[k]{S}\ma{\Delta}-\tfrac14\mat{\Delta}\ma[i]{S}\mat[i]{S}\ma[k]{S}\mat{\Delta}\ma[k]{S}-\\
&\hphantom{a_{ik}\trace(}\tfrac14\mat[i]{S}\ma{\Delta}\mat[i]{S}\ma[k]{S}\mat[k]{S}\ma{\Delta}-\tfrac14\mat[i]{S}\ma{\Delta}\mat[i]{S}\ma[k]{S}\mat{\Delta}\ma[k]{S}).
\end{align*}
Use the identity $\trace\ma{ABCD}=\langle\vect\mat{A},(\mat{D}\otimes\ma{B})\vect\ma{C}\rangle$ \citep{graham1981kronecker} and the notation $\ve[1]{d}=\vect\ma{\Delta}$, $\ve[2]{d}=\vect\mat{\Delta}$ to write
%
\begin{align*}
q_{ik}={}&a_{ik}(-\|\ve[1]{d}\|^2+\tfrac12\langle\ve[1]{d},(\ma[p]{I}\otimes\ma[i]{S}\mat[i]{S})\ve[1]{d}\rangle+\\
&\tfrac12\langle\ve[2]{d},(\ma[i]{S}\otimes\mat[i]{S})\ve[1]{d}\rangle+\tfrac12\langle\ve[1]{d},(\ma[p]{I}\otimes\ma[k]{S}\mat[k]{S})\ve[1]{d}\rangle+\\
&\tfrac12\langle\ve[1]{d},(\mat[k]{S}\otimes\ma[k]{S})\ve[2]{d}\rangle-\tfrac14\langle\ve[1]{d},(\ma[p]{I}\otimes\ma[i]{S}\mat[i]{S}\ma[k]{S}\mat[k]{S})\ve[1]{d}\rangle-\\
&\tfrac14\langle\ve[1]{d},(\mat[k]{S}\otimes\ma[i]{S}\mat[i]{S}\ma[k]{S})\ve[2]{d}\rangle-\\
&\tfrac14\langle\ve[2]{d},(\ma[i]{S}\otimes\mat[i]{S}\ma[k]{S}\mat[k]{S})\ve[1]{d}\rangle-\tfrac14\langle\ve[2]{d}(\ma[i]{S}\mat[k]{S}\otimes\mat[i]{S}\ma[k]{S})\ve[2]{d}\rangle)\\
={}&\langle\ve{d},\ma[ik]{Q}\ve{d}\rangle,
\end{align*}
where $\ma[ik]{Q}$ is given in Table \ref{tab:matrix} and $\ve{d}=[\vet[1]{d}\,\vet[2]{d}]\mtr$.
\begin{table*}
	\begin{align*}
	\ma[ik]{Q}={}&a_{ik}\begin{bmatrix}
	-\ma[np]{I}+\tfrac12\ma[p]{I}\otimes(\ma[i]{S}\mat[i]{S}+\ma[k]{S}\mat[k]{S})-\tfrac14\ma[p]{I}\otimes\ma[i]{S}\mat[i]{S}\ma[k]{S}\mat[k]{S} & \tfrac12\mat[k]{S}\otimes\ma[k]{S}-\tfrac14\mat[k]{S}\otimes\ma[i]{S}\mat[i]{S}\ma[k]{S}\\
	\tfrac12\ma[i]{S}\otimes\mat[i]{S}-\tfrac14\ma[i]{S}\otimes\mat[i]{S}\ma[k]{S}\mat[k]{S} & -\tfrac14\ma[i]{S}\mat[k]{S}\otimes\mat[i]{S}\ma[k]{S}
	\end{bmatrix}
	\end{align*}
	\caption{The matrix $\ma[ik]{Q}$.\label{tab:matrix}}
\end{table*}
%

Furthermore,
\begin{align*}
q_i={}&\langle\Pi_i\ma{\Delta},\ma[i]{\Sigma}\syms\mat[i]{V}\Pi_i\ma{\Delta}+\Pi_i(\ma{\Delta})\ma[i]{\Gamma}\rangle\\
={}&\langle\Pi_i\ma{\Delta},\ma[i]{\Sigma}\syms\mat[i]{V}\Pi_i\ma{\Delta}\rangle+\langle\Pi_i\ma{\Delta},\Pi_i(\ma{\Delta})\ma[i]{\Gamma}\rangle.
\end{align*}
Since $\langle\mat[i]{S}\Pi_i\ma{\Delta},\syms\mat[i]{\Sigma}\Pi_i\ma{\Delta}\rangle=\ma{0}$ by the orthogonality of symmetric and skew-symmetric matrices, we get 
\begin{align*}
q_i={}&\langle\Pi_i\ma{\Delta},\Pi_i(\ma{\Delta})\ma[i]{\Gamma}\rangle\\
={}&\langle\ma{\Delta}-\ma[i]{S}\syms\mat[i]{S}\ma{\Delta},(\ma{\Delta}-\ma[i]{S}\syms\mat[i]{S}\ma{\Delta})\ma[i]{\Gamma}\rangle\\
={}&\trace(\mat{\Delta}\ma{\Delta}-2\syms(\mat[i]{S}\ma{\Delta})\mat[i]{S}\ma{\Delta}+(\syms\mat[i]{S}\ma{\Delta})^2)\ma[i]{\Gamma}\\
={}&\trace(\mat{\Delta}\ma{\Delta}-(\mat[i]{S}\ma{\Delta}+\mat{\Delta}\ma[i]{S})\mat[i]{S}\ma{\Delta}+\\
&\hphantom{\trace(}\tfrac14(\mat[i]{S}\ma{\Delta}+\mat{\Delta}\ma[i]{S})(\mat[i]{S}\ma{\Delta}+\mat{\Delta}\ma[i]{S}))\ma[i]{\Gamma}\\
={}&\trace(\mat{\Delta}\ma{\Delta}-\mat[i]{S}\ma{\Delta}\mat[i]{S}\ma{\Delta}-\mat{\Delta}\ma[i]{S}\mat[i]{S}\ma{\Delta}+\\
&\hphantom{\trace(}\tfrac14(\mat[i]{S}\ma{\Delta}\mat[i]{S}\ma{\Delta}+\mat[i]{S}\ma{\Delta}\mat{\Delta}\ma[i]{S}+\\
&\hphantom{\trace\tfrac14(}\mat{\Delta}\ma[i]{S}\mat[i]{S}\ma{\Delta}+\mat{\Delta}\ma[i]{S}\mat{\Delta}\ma[i]{S}))\ma[i]{\Gamma}\\
={}&\trace(\mat{\Delta}\ma{\Delta}-\tfrac12\mat[i]{S}\ma{\Delta}\mat[i]{S}\ma{\Delta}-\tfrac34\mat{\Delta}\ma[i]{S}\mat[i]{S}\ma{\Delta}+\\
&\hphantom{\trace(}\tfrac14\mat[i]{S}\ma{\Delta}\mat{\Delta}\ma[i]{S})\ma[i]{\Gamma}\\
={}&\trace(\mat{\Delta}\ma{\Delta}\ma[i]{\Gamma}-\tfrac12\ma{\Delta}\mat[i]{S}\ma{\Delta}\ma[i]{\Gamma}\mat[i]{S}-\\
&\hphantom{\trace(}\tfrac34\mat{\Delta}\ma[i]{S}\mat[i]{S}\ma{\Delta}\ma[i]{\Gamma}+\tfrac14\ma{\Delta}\mat{\Delta}\ma[i]{S}\ma[i]{\Gamma}\mat[i]{S})\\
={}&\vet[1]{d}(\ma[i]{\Gamma}\otimes\ma[n]{I})\ve[1]{d}-\tfrac12\vet[2]{d}(\ma[i]{S}\ma[i]{\Gamma}\otimes\mat[i]{S})\ve[1]{d}-\\
&\tfrac34\vet[1]{d}(\ma[i]{\Gamma}\otimes\ma[i]{S}\mat[i]{S})\ve[1]{d}+\tfrac14\vet[2]{d}(\ma[i]{S}\ma[i]{\Gamma}\mat[i]{S}\otimes\ma[p]{I})\ve[2]{d}\\
={}&\langle\ve{d},\ma[i]{Q}\ve{d}\rangle,
\end{align*}
where
\begin{align*}
\ma[i]{Q}={}&\begin{bmatrix}
\ma[i]{\Gamma}\otimes\ma[n]{I}-\tfrac34\ma[i]{\Gamma}\otimes\ma[i]{S}\mat[i]{S} & \ma{0}\\
-\tfrac12\ma[i]{S}\ma[i]{\Gamma}\otimes\mat[i]{S} & \tfrac14\ma[i]{S}\ma[i]{\Gamma}\mat[i]{S}\otimes\ma[p]{I}
\end{bmatrix}.
\end{align*}

There is a constant permutation matrix $\ma{K}\in\mathsf{O}(np)$ such that $\vect\mat{\Delta}=\ma{K}\vect\ma{\Delta}$ for all $\vect\ma{\Delta}\in\R^{np}$ \cite{graham1981kronecker}. Hence
\begin{align*}
\ve{d}=\begin{bmatrix}
\vect\ma{\Delta}\\
\vect\mat{\Delta}\\
\end{bmatrix}=\begin{bmatrix}
\ma[np]{I}\\
\ma{K}
\end{bmatrix}
\vect\ma{\Delta}=\begin{bmatrix}
\ma[np]{I}\\
\ma{K}
\end{bmatrix}
\ve[1]{d}.
\end{align*}
The quadratic form $q$ satisfies
\begin{align*}
q&=\sum_{i\in\V}\langle\ve{d},\ma[i]{Q}\ve{d}\rangle+2\sum_{e\in\E}\langle\ve{d},\ma[ik]{Q}\ve{d}\rangle\\
&=\Bigl\langle\ve{d},\Bigl(\sum_{i\in\V}\ma[i]{Q}+\sum_{k\in\Ni}\ma[ik]{Q}\Bigr)\ve{d}\Bigr\rangle\\
&=\left\langle\begin{bmatrix}
\ma[np]{I}\\
\ma{K}
\end{bmatrix}\ve[1]{d},\Bigl(\sum_{i\in\V}\ma[i]{Q}+\sum_{k\in\Ni}\ma[ik]{Q}\Bigr)\begin{bmatrix}
\ma[np]{I}\\
\ma{K}
\end{bmatrix}\ma[1]{d}\right\rangle\\
&=\Bigl\langle\ve[1]{d},\begin{bmatrix}
\ma[np]{I} & \!\mat{K} \end{bmatrix}\Bigl(\sum_{i\in\V}\ma[i]{Q}+\sum_{k\in\Ni}\ma[ik]{Q}\Bigr)\begin{bmatrix}
\ma[np]{I}\\
\ma{K}
\end{bmatrix}\ve[1]{d}\Bigr\rangle\\
&=\langle\ve[1]{d},\ma{M}\ve[1]{d}\rangle,
\end{align*}
where
\begin{align*}
\ma{M}=\syms\begin{bmatrix}
\ma[np]{I} & \!\mat{K} \end{bmatrix}\ma{Q}\begin{bmatrix}
\ma[np]{I}\\
\ma{K}
\end{bmatrix},\!\quad\ma{Q}=\sum_{i\in\V}\ma[i]{Q}+\sum_{k\in\Ni}\ma[ik]{Q}.
\end{align*}
%


%
%

\subsection{Upper bound of the smallest eigenvalue}

\label{sec:upper}

We wish to show that $q$ assumes negative values for some $\ma{\Delta}\in\R^{n\times p}$ at all equilibria $\ma{S}\notin\mathcal{C}$. This excludes any such equilibria from being a local minimizer of the potential function $V$ given by \eqref{eq:VSt}. If $\trace\ma{M}$ is negative, then $\ma{M}$ has at least one negative eigenvalue. Calculate
\begin{align*}
\trace\ma{M}&=
\trace\syms(\ma{Q})\begin{bmatrix}
\ma[np]{I} & \mat{K}\\
\ma{K} & \ma[np]{I}
\end{bmatrix}
=\trace\ma{A}+2\trace\ma{B}\ma{K}+\trace\ma{C},
\end{align*}
where $\ma{A}$, $\ma{B}$, and $\ma{C}$ denote the three blocks of $\syms\ma{Q}$. Let us calculate each of the three terms in $\trace\ma{M}$ separately, starting with $\ma{A}$ and $\ma{C}$,
\begin{align*}
\trace\ma{A}={}&\sum_{i\in\V}\trace(\ma[i]{\Gamma}\otimes\ma[n]{I}-\tfrac34\ma[i]{\Gamma}\otimes\ma[i]{S}\mat[i]{S})+\\
&\sum_{k\in\Ni}\!\!a_{ik}\trace(-\ma[np]{I}+\tfrac12\ma[p]{I}\otimes(\ma[i]{S}\mat[i]{S}+\ma[k]{S}\mat[k]{S})-\\
&\hphantom{\sum_{k\in\Ni}a_{ik}\trace(}\,\!\tfrac14\ma[p]{I}\otimes\ma[i]{S}\mat[i]{S}\ma[k]{S}\mat[k]{S})\\
={}&\sum_{i\in\V}n\trace\ma[i]{\Gamma}-\tfrac{3p}4\trace\ma[i]{\Gamma}+\\
&\sum_{k\in\Ni}a_{ik}(-np+p^2-\tfrac{p}4\|\mat[k]{S}\ma[i]{S}\|^2),
\end{align*}
where we utilize that
\begin{align*} \trace\ma{X}\otimes\ma{Y}&=\trace\ma{X}\trace\ma{Y},\\ \trace\ma{S}\mat{S}&=\trace\mat{S}\ma{S}=\trace\ma[p]{I}=p,\\
\trace\ma{Z}\mat{Z}\ma{W}\mat{W}&=\trace(\mat{Z}\ma{W})\mtr(\mat{Z}\ma{W})=\|\mat{Z}\ma{W}\|^2,
\end{align*}
for any $\ma{X},\ma{Y}\in\R^{n\times n}$, $\ma{S}\in\St$, and $\ma{Z}\in\R^{n\times p},\ma{W}\in\R^{n\times q}$. Continuing,
\begin{align*}
\trace\ma{A}={}&\!\sum_{i\in\V}\!\left(n-\tfrac{3p}4\right)\trace\ma[i]{\Gamma}-\!\sum_{k\in\Ni}\!a_{ik}((n-p)p+\tfrac{p}4\|\mat[k]{S}\ma[i]{S}\|^2)\\
={}&\sum_{i\in\V}\sum_{k\in\Ni}\!\!\left(n-\tfrac{3p}4\right)\left\langle a_{ik}\ma[k]{S},\ma[i]{S}\right\rangle-\\
&\hspace{1.3cm}a_{ik}((n-p)p-\tfrac{p}4\|\mat[k]{S}\ma[i]{S}\|^2)\\
={}&2\sum_{e\in\E}a_{ik}\!\left(\left(n-\tfrac{3p}4\right)\left\langle\ma[k]{S},\ma[i]{S}\right\rangle-\right.\\
&\hspace{1.45cm}(n-p)p-\tfrac{p}4\|\mat[k]{S}\ma[i]{S}\|^2\bigr),\\
\trace\ma{C}={}&\sum_{i\in\V}\tfrac14\trace(\ma[i]{S}\ma[i]{\Gamma}\mat[i]{S}\otimes\ma[p]{I})-\!\sum_{k\in\Ni}\tfrac{a_{ik}}4\trace(\ma[i]{S}\mat[k]{S}\otimes\mat[i]{S}\ma[k]{S})\\
={}&\sum_{i\in\V}\tfrac{p}4\trace(\ma[i]{\Gamma})-\sum_{k\in\Ni}\tfrac{a_{ik}}4\trace(\mat[k]{S}\ma[i]{S})^2\\
={}&\sum_{i\in\V}\sum_{k\in\Ni}\tfrac{p}4\left\langle a_{ik}\ma[k]{S},\ma[i]{S}\right\rangle-\tfrac{a_{ik}}4\langle\ma[k]{S},\ma[i]{S}\rangle^2\\
={}&2\sum_{e\in\E}a_{ik}\bigl(\tfrac{p}4\left\langle\ma[k]{S},\ma[i]{S}\right\rangle-\tfrac{1}4\langle\ma[k]{S},\ma[i]{S}\rangle^2\bigr).
\end{align*}
Note that 
\begin{align*}
\ma{B}={}&\sum_{i\in\V}-\tfrac14\ma[i]{\Gamma}\mat[i]{S}\otimes\ma[i]{S}+\\
&\hspace{0.5cm}\sum_{k\in\Ni}a_{ik}(\tfrac14\mat[k]{S}\otimes\ma[k]{S}-\tfrac18\mat[k]{S}\otimes\ma[i]{S}\mat[i]{S}\ma[k]{S}+\\
&\hspace{1.8cm}\tfrac14\mat[i]{S}\otimes\ma[i]{S}-\tfrac18\mat[i]{S}\otimes\ma[k]{S}\mat[k]{S}\ma[i]{S}).
\end{align*}
To calculate $\trace\ma{B}\ma{K}$, we utilize that $\ma{K}=\sum_{a=1}^n\sum_{b=1}^p\ma[ab]{E}\otimes\ma[ba]{E}$, where the elemental matrix $\ma[ab]{E}\in\R^{n\times p}$ is given by $\ma[ab]{E}=\ve[a]{e}\otimes\ve[b]{e}$ for all $a\in\{1,\ldots,n\}$, $b\in\{1,\ldots,p\}$ \citep{graham1981kronecker}:
\begin{align*}
\trace\ma{B}\ma{K}={}&\sum_{i\in\V}-\tfrac14\trace(\ma[i]{\Gamma}\mat[i]{S}\otimes\ma[i]{S})\ma{K}+\\
&\sum_{k\in\Ni}a_{ik}\trace(\tfrac14\mat[k]{S}\otimes\ma[k]{S}-\tfrac18\mat[k]{S}\otimes\ma[i]{S}\mat[i]{S}\ma[k]{S})\ma{K}+\\
&\sum_{k\in\Ni}a_{ik}\trace(\tfrac14\mat[i]{S}\otimes\ma[i]{S}-\tfrac18\ma[i]{S}\otimes\ma[k]{S}\mat[k]{S}\ma[i]{S})\ma{K}\\
={}&\sum_{i\in\V}\sum_{a,b}-\tfrac14\trace(\ma[i]{\Gamma}\mat[i]{S}\ma[ab]{E}\otimes\ma[i]{S}\ma[ba]{E})+\\
&\sum_{k\in\Ni}a_{ik}\sum_{a,b}\trace\bigl(\tfrac14\mat[k]{S}\ma[ab]{E}\otimes\ma[k]{S}\ma[ba]{E}-\\
&\hspace{2.1cm}\tfrac18\mat[k]{S}\ma[ab]{E}\otimes\mat[i]{S}\mat[i]{S}\ma[k]{S}\ma[ba]{E}\bigr)+\\
&\sum_{k\in\Ni}a_{ik}\sum_{a,b}\trace\bigl(\tfrac14\mat[i]{S}\ma[ab]{E}\otimes\ma[i]{S}\ma[ba]{E}-\\
&\hspace{2.1cm}\tfrac18\mat[i]{S}\ma[ab]{E}\otimes\ma[k]{S}\mat[k]{S}\ma[i]{S}\ma[ba]{E}\bigr),
\end{align*}
where we use the mixed-product property of Kronecker products, $(\ma{X}\otimes\ma{Y})(\ma{Z}\otimes\ma{W})=(\ma{X}\ma{Z})\otimes(\ma{Y}\ma{W})$, which holds for any matrices $\ma{X},\ma{Y},\ma{Z},\ma{W}$ such that $\ma{X}\ma{Z}$ and $\ma{Y}\ma{W}$ are well-defined. Continuing,
\begin{align*}
\trace\ma{B}\ma{K}={}&\sum_{i\in\V}\sum_{a,b}-\tfrac14\trace(\ma[i]{\Gamma}\mat[i]{S}\ma[ab]{E})\trace(\ma[i]{S}\ma[ba]{E})+\\
&\sum_{k\in\Ni}a_{ik}\sum_{a,b}\tfrac14\trace(\mat[k]{S}\ma[ab]{E})\trace(\ma[k]{S}\ma[ba]{E})-\\
&\hspace{2.1cm}\tfrac18\trace(\mat[k]{S}\ma[ab]{E})\trace(\ma[i]{S}\mat[i]{S}\ma[k]{S}\ma[ba]{E})+\\
&\sum_{k\in\Ni}a_{ik}\sum_{a,b}\tfrac14\trace(\mat[i]{S}\ma[ab]{E})\trace(\ma[i]{S}\ma[ba]{E})-\\
&\hspace{2.1cm}\tfrac18\trace(\mat[i]{S}\ma[ab]{E})\trace(\ma[k]{S}\mat[k]{S}\ma[i]{S}\ma[ba]{E})\\
={}&\sum_{i\in\V}\sum_{a,b}-\tfrac14\trace\Bigl(\sum_{k\in\Ni}a_{ik}\mat[k]{S}\ma[ab]{E}\Bigr)\trace(\ma[i]{S}\ma[ba]{E})+\\
&\sum_{k\in\Ni}a_{ik}\sum_{a,b}\tfrac14(\ma[k]{S})_{ab}(\ma[k]{S})_{ab}-\\
&\hspace{2.1cm}\tfrac18(\ma[k]{S})_{ab}(\ma[i]{S}\mat[i]{S}\ma[k]{S})_{ab}+\\
&\sum_{k\in\Ni}a_{ik}\sum_{a,b}\tfrac14(\ma[i]{S})_{ab}(\ma[i]{S})_{ab}-\\
&\hspace{2.1cm}\tfrac18(\ma[i]{S})_{ab}(\ma[k]{S}\mat[k]{S}\ma[i]{S})_{ab}\\
={}&\sum_{i\in\V}\sum_{a,b}\sum_{k\in\Ni}-\tfrac{a_{ik}}4\trace(\ma[k]{S}\ma[ba]{E})\trace(\ma[i]{S}\ma[ba]{E})+\\
&\sum_{k\in\Ni}a_{ik}\bigl(\tfrac14\|\ma[k]{S}\|^2-\tfrac18\langle\ma[k]{S},\ma[i]{S}\mat[i]{S}\ma[k]{S}\rangle\bigr)+\\
&\sum_{k\in\Ni}a_{ik}\bigl(\tfrac14\|\ma[i]{S}\|^2-\tfrac18\langle\ma[i]{S},\ma[k]{S}\mat[k]{S}\ma[i]{S}\rangle\bigr),
\end{align*}
where we utilize that
\begin{align*}
\sum_{a,b}(\ma[ab]{X})^2=\|\ma{X}\|^2, \,\, \sum_{a,b}\ma[ab]{X}\ma[ab]{Y}=\langle\ma{X},\ma{Y}\rangle
\end{align*} 
for all $\ma{X},\ma{Y}\in\R^{n\times m}$. Finally,
\begin{align*}
\trace\ma{B}\ma{K}={}&\sum_{i\in\V}\sum_{k\in\Ni}-\tfrac{a_{ik}}4\langle\ma[k]{S},\ma[i]{S}\rangle+a_{ik}\bigl(\tfrac{p}2-\tfrac14\|\mat[i]{S}\ma[k]{S}\|^2\bigr)\\
={}&2\sum_{e\in\E}a_{ik}\bigl(-\tfrac14\langle\ma[k]{S},\ma[i]{S}\rangle+\tfrac{p}2-\tfrac14\|\mat[i]{S}\ma[k]{S}\|^2\bigr).
\end{align*}

Adding up all four terms gives
\begin{align}
\tfrac12\trace\ma{M}={}&\tfrac{1}{2}\trace\ma{A}+\trace\ma{B}\ma{K}+\tfrac12\trace\ma{C}\nonumber\\
={}&\sum_{e\in\E}a_{ik}\bigl(\left(n-\tfrac{3p}4\right)\left\langle\ma[k]{S},\ma[i]{S}\right\rangle-(n-p)p\nonumber\\
&\hspace{1.2cm}-\tfrac{p}4\|\mat[k]{S}\ma[i]{S}\|^2-\tfrac12\langle\ma[k]{S},\ma[i]{S}\rangle+p-\nonumber\\
&\hspace{1.3cm}\tfrac12\|\mat[i]{S}\ma[k]{S}\|^2+\tfrac{p}4\left\langle\ma[k]{S},\ma[i]{S}\right\rangle-\tfrac{1}4\langle\ma[k]{S},\ma[i]{S}\rangle^2\bigr)\nonumber\\
={}&\sum_{e\in\E}a_{ik}\bigl(\left(n-\tfrac{p+1}2\right)\left\langle\ma[k]{S},\ma[i]{S}\right\rangle-\tfrac{p+2}4\|\mat[k]{S}\ma[i]{S}\|^2-\nonumber\\
&\hphantom{\sum_{e\in\E}a_{ik}\bigl(}\,\,\,\tfrac14\langle\ma[k]{S},\ma[i]{S}\rangle^2+(1-n+p)p\bigr).\label{eq:trM}
\end{align}

Equation \eqref{eq:trM} is the desired expression for $\trace\ma{M}$. In the next section we will study how it varies over $\St^N$. To verify that no miscalculations were made, note that at a consensus, where $\mat[k]{S}\ma[i]{S}=\ma[p]{I}$, we get 
\begin{align*}
\tfrac12\trace\ma{M}|_{\mathcal{C}}&=\sum_{e\in\E}a_{ik}\bigl(n-\tfrac{p+1}{2}-\tfrac{p+2+p}{4}+1-n+p\bigr)p=0
\end{align*}
This is expected since $\mathcal{C}$ is invariant under any tangent vector that belongs to its tangent space, $\ma[i]{\Delta}|_{\mathcal{C}}=(\Pi_1\ma{\Delta})_{i=1}^N\in\ts[\St^N]{\mathcal{C}}$, and $V$ is constant over $\mathcal{C}$. Also note that \eqref{eq:trM} is consistent with the corresponding expression in \cite{markdahl2018tac} for the special case of $\Sn=\mathsf{St}(1,n+1)$.

\subsection{Nonlinear programming problem}
\label{sec:nlp}

It remains to show that $\trace\ma{M}$ given by \eqref{eq:trM} is strictly negative for each equilibrium configuration $\ma{S}\notin\mathcal{C}$. To that end, we could consider the problem of maximizing $\trace\ma{M}$ over all configurations $\ma{S}\notin\mathcal{C}$ which satisfy the equations \eqref{eq:eq} that characterize an equilibrium set. However, that problem seems difficult to solve. Instead, we make use of the following inequality
\begin{align}
\begin{split}\label{eq:opt}
\tfrac12\trace\ma{M}\leq{}& |\E|\max_{e\in\E}a_{ik}\max_{\ma{X},\ma{Y}}f(\ma{X},\ma{Y}),\\
f(\ma{X},\ma{Y})={}&\left(n-\tfrac{p+1}2\right)\left\langle\ma{X},\ma{Y}\right\rangle-\tfrac{p+2}4\|\mat{X}\ma{Y}\|^2-\\
&\tfrac14\langle\ma{X},\ma{Y}\rangle^2+(1-n+p)p,
\end{split}
\end{align}
where $f:\St\times\St\rightarrow\R$. If we can show that the upper bound on $\trace\ma{M}$ is negative for all $\ma{X}\neq\ma{Y}$, then we are done. Note that the inequality is sharp in the case of two agents and that $f(\ma{X},\ma{X})=0$ since this corresponds to consensus in a system of two agents.

Denote $\ma{Z}=\mat{X}\ma{Y}$. It is clear that $\trace\ma{Z}\in[-p,p]$ since
\begin{align*}
|\trace\ma{Z}|\leq\Bigl|\sum_{i=1}^p\lambda_i\Bigr|\leq p\|\ma{Z}\|_2\leq p\|\ma{X}\|_2\|\ma{Y}\|_2=p.
\end{align*}
Consider a relaxation of \eqref{eq:opt} where $\ma{Z}\in\R^{n\times p}$ subject to $\trace\ma{Z}\in[-p,p]$. Let $\widebar{f}:\R^{p\times p}\rightarrow\R$ denote the extension of $f$ given by
\begin{align*}
\widebar{f}(\ma{Z})={}&(n-\tfrac{p+1}{2})\trace\ma{Z}-\tfrac{p+2}{4}\|\ma{Z}\|^2-\tfrac14(\trace\ma{Z})^2\\
&+(1-n+p)p.
\end{align*}
Note that $\bar{f}$ being negative for all $\ma{Z}\in\R^{n\times p}$ with $\trace\ma{Z}\in[-p,p]$, implies that $f(\ma{X},\ma{Y})$ is negative for all $\ma{X},\ma{Y}\in\St$. To simplify $\widebar{f}(\ma{Z})$, first observe that
\begin{align*}
\|\ve{Z}\|^2&\geq\sum_{i=1}^p|\lambda_i|^2=\|[\lambda_i]\|^2\geq\tfrac1p\|[\lambda_i]\|_1^2\\
&=\tfrac1p\Bigl(\sum_{i=1}^p|\lambda_i|\Bigr)^2\geq\tfrac1p\Bigl(\sum_{i=1}^p|\re\lambda_i|\Bigl)^2\\
&\geq\tfrac1p\Bigl(\sum_{i=1}^p\re\lambda_i\Bigr)^2=\tfrac{1}{p}(\trace\ma{Z})^2,
\end{align*}
where Schur's inequality relates the Frobenius norm of $\ma{Z}$ to its eigenvalues $\lambda_1,\ldots,\lambda_p$ \citep{horn2012matrix}. Use the above inequality to write
\begin{align*}
\widebar{f}(\ma{Z})\leq(n-\tfrac{p+1}{2})\trace\ma{Z}-\tfrac{p+1}{2p}(\trace\ma{Z})^2+(1-n+p)p.
\end{align*}
Note that this upper bound on $\bar{f}(\ma{Z})$ is quadratic in $\trace\ma{Z}$. The maximum of the quadratic is located at 
\begin{align*}
\trace\ma{Z}=\tfrac{(2n-p-1)p}{2(p+1)}.
\end{align*}
Assume that the maximum of the parabola is larger than $p$, \ie $\trace\ma{Z}\geq p$. Simplifying this inequality we find that $p\leq\tfrac23n-1$. Since the bound is a concave quadratic polynomial, its maximum value for $\trace\ma{Z}\in[-p,p]$ is obtained at the feasible point that is closest to the optimal point, \ie at $\trace\ma{Z}=p$ where the bound equals $0$. The value $\trace\ma{Z}=p$ can only be achieved when $\ma{Y}=\ma{X}$.

\end{document}